\documentclass{amsart}

\usepackage{enumerate,url,amssymb}

\newtheorem{theorem}{Theorem}[section]
\newtheorem{lemma}[theorem]{Lemma}
\newtheorem{proposition}[theorem]{Proposition}

\theoremstyle{definition}

\newtheorem{example}[theorem]{Example}
\newtheorem{corollary}[theorem]{Corollary}

\theoremstyle{remark}

\numberwithin{equation}{section}

\newcommand{\abs}[1]{\lvert#1\rvert}
\newcommand{\dtext}{\textnormal d}

\DeclareMathOperator{\re}{Re}
\DeclareMathOperator{\im}{Im}
\DeclareMathOperator{\id}{id}
\DeclareMathOperator{\loc}{loc}
\def\Xint#1{\mathchoice
{\XXint\displaystyle\textstyle{#1}}%
{\XXint\textstyle\scriptstyle{#1}}%
{\XXint\scriptstyle\scriptscriptstyle{#1}}%
{\XXint\scriptscriptstyle\scriptscriptstyle{#1}}%
\!\int}
\def\XXint#1#2#3{{\setbox0=\hbox{$#1{#2#3}{\int}$}
\vcenter{\hbox{$#2#3$}}\kern-.5\wd0}}

\def\dashint{\Xint-}
\def\le{\leqslant}
\def\ge{\geqslant}

\begin{document}

\title{On Injectivity of Quasiregular Mappings}

\author{Tadeusz Iwaniec}
\address{Department of Mathematics, Syracuse University, Syracuse,
NY 13244, USA}
\email{tiwaniec@syr.edu}
\thanks{Iwaniec was supported  by the NSF grant DMS-0800416.}

\author{Leonid V. Kovalev}
\address{Department of Mathematics, Syracuse University, Syracuse,
NY 13244, USA}
\email{lvkovale@syr.edu}
\thanks{Kovalev was supported by the NSF grant DMS-0700549.}

\author{Jani Onninen}
\address{Department of Mathematics, Syracuse University, Syracuse,
NY 13244, USA}
\email{jkonnine@syr.edu}
\thanks{Onninen was supported by the NSF grant  DMS-0701059.}


\date{September 28, 2008}


\begin{abstract}
We give sufficient conditions for a planar quasiregular mapping to be injective in terms of the range of the differential matrix.
\end{abstract}

\maketitle

\section{Introduction}

A holomorphic function $f\colon \Omega\to\mathbb C$ of one complex variable is a local homeomorphism if and only if $f'\ne 0$
in $\Omega$. If, in addition, $\Omega$ is convex and $\re f'\geqslant 0$, then $f$ is either constant or  injective in $\Omega$~\cite{Al}.
In this paper we establish the analogues of these well-known facts for quasiregular mappings $f\colon \Omega\to\mathbb C$. By definition
 $f\in W_{\rm loc}^{1,2}(\Omega)$ is quasiregular is there exists a constant $k<1$ such that  $\abs{f_{\bar z}}\leqslant k\abs{f_{z}}$ a.e. in $\Omega$.  Such a mapping can be called $K$-quasiregular with $K=(1+k)/(1-k)$, or $K$-quasiconformal if it is also injective.

\begin{theorem}\label{T2}
 Let $\Omega \subset \mathbb C$ be a domain. If  $f :\Omega  \rightarrow \mathbb C$ is a non-constant  quasiregular mapping and
 $\re f_z \geqslant 0$ almost everywhere, then $f$ is  a local homeomorphism.
 \end{theorem}

The proof is based on the celebrated theorem of Poincar\'e-Bendixson \cite{Be}  and its extension by Brouwer \cite{Br}, about local structure of integral curves of a continuous planar vector field near its critical point. Example~\ref{branchexample} will show that the assumption $\re f_z\geqslant 0$ cannot be replaced with $\abs{\arg f_z}\leqslant \pi/2+\epsilon$, for any $\epsilon>0$.  To ensure that $f$ is injective in a convex domain $\Omega \subsetneq \mathbb C$ we must restrict the range
of $f_z$  even further.

\begin{theorem}\label{T15}
Let $\Omega \subset \mathbb C$ be a convex domain. If $f\colon \Omega  \rightarrow \mathbb C$ is a non-constant  quasiregular mapping and $\re f_z = 0$ almost everywhere, then $f$ is  a homeomorphism.
\end{theorem}

This theorem admits the following reformulation (see section~\ref{PT15}): if $\psi$ is a differentiable
real-valued function on a convex domain $\Omega\subset \mathbb R^2$ and the gradient mapping $\nabla \psi\colon\Omega\to\mathbb R^2$ is quasiregular, then $\nabla \psi$ is either injective or constant. This is no longer true in dimensions $n\ge 3$, as is demonstrated by  Example~\ref{gradexample}.
Also, the assumption $\re f_z=0$ in Theorem~\ref{T15} cannot be replaced with $\abs{\arg f_z}<\epsilon$, for any $\epsilon>0$, by
Example~\ref{noninjexample}. However, the situation is different when $\Omega=\mathbb C$. This can be expected since by Picard's
theorem an entire function whose derivative omits two values is linear, and therefore is either constant or injective. For quasiregular
mappings we have the following

\begin{theorem}\label{T1}
If $f\colon\mathbb C \rightarrow \mathbb C$ is a non-constant  quasiregular mapping and $\re  f_z \geqslant 0$ almost everywhere, then $f$ is  a homeomorphism.
\end{theorem}

The sharpness of the assumption is demonstrated by Example~\ref{branchexample}. As a corollary of Theorem~\ref{T2} and~\ref{T1} we obtain a converse to the following theorem \cite[Theorem 6.3.1]{AIMbook}:
\begin{theorem}
If $f\in W^{1,2}_{\loc} (\mathbb C)$ is a homeomorphic solution to the reduced Beltrami equation
\begin{equation}\label{RBE}
f_{\bar z} = \lambda (z)  \re f_z \, , \hskip1cm  |\lambda (z)| \leqslant k <1 \,
\end{equation}
then $\re f_z$ does not change sign.
\end{theorem}
\begin{corollary}
If $f\in W^{1,2}_{\loc} (\mathbb C)$ is a  solution of (\ref{RBE}) such that  $\re f_z$ does not change sign, then $f$ is a homeomorphism.
\end{corollary}
Let us emphasize that the notion of quasiregularity is invariant under affine change of variables. Accordingly, the assumption $\re f_z \geqslant 0$ in Theorems \ref{T2} and \ref{T1} can be replaced by somewhat geometrically pleasing condition on the differential matrix $Df(z)$, which is also invariant under affine change of variables. Let $\mathbb R^{2\times 2}$ denote the $4$-dimensional linear space of $2 \times 2$ matrices equipped with the inner product $\langle X, Y \rangle = \textnormal{Tr}\, (X \, Y^\ast)$, for $X,Y \in \mathbb R^{2 \times 2}$. Each nonzero matrix $N\in \mathbb R^{2 \times 2}$ gives rise to a $3$-dimensional subspace perpendicular to $N$,
$$\mathbb H_N=\big\{X \colon \langle X, N \rangle =0 \big\} \subset \mathbb R^{2 \times 2}$$
There are three types of such subspaces:
\begin{itemize}
\item{$\mathbb H_N$ is  said to be {\it positive  subspace} if $\det N >0$  }
\item{$\mathbb H_N$ is  said to be {\it negative subspace} if $\det N <0$  }
\item{$\mathbb H_N$ is said to be {\it singular} if $\det N =0$}
\end{itemize}
A $3$-dimensional subspace $\mathbb H_N$ splits $\mathbb R^{2 \times 2}$ into two half-spaces. If $\det N >0$
we call them {\it positive half-spaces} of $\mathbb R^{2\times 2}$. If $\det N <0$ we call them {\it negative half-spaces}.  If $\det N=0$ the corresponding half-spaces are called {\it singular half-spaces}. Now, upon affine change of variable the assumption $\re f_z \geqslant 0 $ tells us that the essential range of $Df$ lies in a positive half-space. Precisely, it means that there is a constant matrix $N\in \mathbb R^{2 \times 2}$ of positive determinant such that
\begin{equation}\label{IEq1}
\langle Df(z), N \rangle \geqslant 0 \hskip2cm \textnormal{ a.e. in }\; \Omega
\end{equation}
As a matter of fact this amounts to saying that the homotopy between $f$ and the $\mathbb R$-linear map $L \colon \mathbb R^2 \rightarrow \mathbb R^2$, with $DL=N$
$$f^t=(1-t)f + tL\, , \hskip2cm 0 \leqslant t \leqslant 1$$
keeps the distortion function of $f^t$ decreasing as $t$ increases from $0$ to $1$. For example, if $L=\id : \mathbb R^2 \rightarrow \mathbb R^2$ then condition (\ref{IEq1}) reads as $\re f_z \geqslant 0$, so $f^t (z)= (1-t)f(z) + tz $ and
\[\left|\frac{f^t_{\bar z  }}{f^t_z   }\right|^2 = \frac{\left| f_{\bar z}\right|^2}{ \left(\re f_z +\frac{t}{1-t}\right)^2 + \left(\im f_z\right)^2} \searrow 0 \hskip1cm \textnormal{ as } \; t \nearrow 1
\]
The limit map $f^1(z)=z$ is a homeomorphism. Recall Hurwitz-type theorems for quasiregular mappings,
see~\cite[II 5.3]{LV} and~\cite[Lemma 3]{Mi}.
\begin{theorem}\label{HuT}
If $f_j : \Omega \rightarrow \mathbb C$ is a sequence of [locally] $K$-quasiconformal mappings which converges  uniformly on compact sets to $f: \Omega  \rightarrow \mathbb C$, then $f$ is either constant or [locally] $K$-quasiconformal.
\end{theorem}
Now, heuristically, by virtue of Theorem \ref{HuT}  it should come by no surprise that Condition (\ref{IEq1}) yields local injectivity of $f$. However, our proof still requires the Poincar\'e-Bendixson analysis  of the integral curves of the vector fields $f^t$, $0<t\leqslant 1$. In view of these observations Theorem~\ref{T2} is a statement on differential inclusions; let
$$\mathbb U\, (K,N)= \big\{X\in \mathbb R^{2 \times 2} \colon |X|^2 \leqslant K \det X \, , \; \; \langle X, N \rangle \geqslant 0 \big\}$$
for some $K \geqslant 1$ and $N \in \mathbb R^{2 \times 2}$ with positive determinant.

 {\it Every nonconstant solution to the differential inclusion $$Df(z)\in \mathbb U\, (K,N)  \quad \textnormal{ a.e. in } \; \Omega, \quad f \in W^{1,2}_{\loc}(\Omega)$$
is a local homeomorphism.}

We refer the reader to \cite{KMS} for a survey on differential inclusions.

\section{Proof of Theorem \ref{T2}}
Let $f$ and $\Omega$ be as in the statement of Theorem \ref{T2}. By virtue of Theorem \ref{HuT}, it suffices
to prove that $f^{\lambda}(z):=f(z)+\lambda z$ is a local homeomorphism for $\lambda>0$. To simplify notation, we write $f$ instead of
$f^{\lambda}$, keeping in mind that $\re f_z\geqslant \lambda>0$ a.e. in $\Omega$.
The local index of $f$ at $z_0\in\Omega$ is an integer defined by the rule
$$n_f(z_0)= \frac{1}{2 \pi} \, \underset{0 \leqslant \theta \leqslant 2 \pi}{\Delta} \arg\left[f(z_0 + re^{i\theta})-f(z_0)\right] $$
where the increment of the argument of $f$ does not depend on the choice of radius $r$, provided $r$ is sufficiently small. It is a general topological fact that $f$ is locally injective if and only if $n_f(z_0)=1$ for every $z_0 \in \Omega$. Since nonconstant   quasiregular mappings are orientation-preserving~\cite[Theorem I.4.5]{Rib}, we have $n_f(z_0) \ge 1$. It remains to show that $n_f(z_0) \leqslant 1$ for every $z_0 \in \Omega$.  It involves no loss of generality in assuming that $z_0=0$, $f(z_0)=0$, $\Omega$
 is the open unit disk, and that $f(z) \not= 0$ for $0<|z| \leqslant 1$. Let $\Omega_\circ $ denote the punctured unit disk, $\Omega_\circ =\{z \colon 0<\abs{z}<1\}$. We shall consider the integral curves of $f$ in $\Omega_\circ$; that is, solutions of the differential  equation
 \begin{equation}
 \frac{\dtext z}{\dtext t} = f(z) \colon 0<|z(t)|<1 \; \; \textnormal{ for }\; a<t<b
  \end{equation}
By virtue of Peano's Existence Theorem, through every point $z_0 \in \Omega_\circ$ there passes an integral curve, though uniqueness is not guaranteed. However, $z=z(t)$ can be extended (as a solution) over a maximal interval of existence, say $a_- < t <b_+$, $-\infty \leqslant a_- <b_+ \leqslant \infty$. Moreover $z(t)$ tends to $\partial \Omega_\circ=\{0\} \cup \mathbb S^1$ as $t \rightarrow a_-$ and $t \rightarrow b_+$. The extension of $z(t)$ need not be unique and the maximal interval of existence depends on the extension. Clearly $z=z(t)$ is of class $C^{1, \alpha} (a_-, b_+)$, $0<\alpha <1$. Since $f(0)=0$, it is possible in general that there is an injective solution $z=z(t)\in \Omega_\circ$ for $0\leqslant t \le \delta$
such that $\lim\limits_{t \rightarrow \delta}z(t)=z(0)$. Then $\gamma=\{z(t) \colon 0 \leqslant t < \delta\}$ is a rectifiable Jordan curve in $\Omega_\circ$.  However under our assumption such curves are not present. Indeed, suppose such $\gamma$ exists. To   reach a contradiction, we let $\mathbb U$ denote the bounded component of $\mathbb C \setminus \gamma$; it is a simply connected region in $\Omega$. We integrate $f_z$ over $\Omega$ by using Stokes' theorem
\begin{eqnarray}
\iint_\mathbb U f_z\, \dtext x \, \dtext y = \frac{1}{2i} \iint_{\mathbb U} \dtext \left(f \, \dtext \overline{z}\right) = \frac{1}{2i} \int_\gamma f\, \dtext \bar z = \frac{1}{2i} \int_0^\delta \left|f\big(z(t)\big)\right|^2 \, \dtext t
\end{eqnarray}
This shows that
\begin{equation}
\iint_\mathbb U \left(\re f_z \right)\, \dtext x \, \dtext y =0
\end{equation}
which is impossible since $\re f_z >0$ almost everywhere.

 Next we shall rule out the integral curves $\gamma = \{z(t)\, ; \; \; a_- <t<b_+\}$ such that $\lim\limits_{t \rightarrow a_-} z(t)= \lim\limits_{t \rightarrow b_+} z(t)=0$. Call such curves {\it elliptic loops} in $\Omega_\circ$.
 According to the celebrated Poincar\'e-Bendixson-Brouwer Theory  \cite{Br} such curves are present in every elliptic sector of $\Omega$. We shall not give a definition of an elliptic sector here as the need will not arise. The interested reader is referred to \cite{Be,Br, Habook} for the definition and thorough discussion of sectors. The proof of nonexistence of elliptic loops is much the same as above. Adding the  point $0$ to $\gamma$ we obtain a Jordan curve, closure of $\gamma$ in $\Omega_\circ$. Let $\mathbb U$ denote the bounded component of $\mathbb C \setminus \overline{\gamma}$. To avoid delicate questions of rectifiability of $\gamma$ we remove from $\mathbb U$ a small disk $\mathbb D_\epsilon=\{z \, ; \; \; |z| \leqslant \epsilon\} \subset \mathbb U$. Then we integrate as before
 \begin{equation}
 \iint_{\mathbb U \setminus \mathbb D_\epsilon} f_z\, \dtext x \, \dtext y = \frac{1}{2i} \int_{\partial (\mathbb U \setminus \mathbb D_\epsilon)} f \, \dtext \bar z = \frac{1}{2i} \int_{\gamma_1} f \, \dtext \bar z +  \frac{1}{2i} \int_{\gamma_2} f \, \dtext \bar z
 \end{equation}
where $\gamma_1 = \partial \mathbb U \setminus \mathbb D_\epsilon$ and $\gamma_2 = \mathbb U \cap \partial  \mathbb D_\epsilon$. As before, the real part of the first integral term vanishes. The second term can be made as small as we wish. Indeed, we have
\begin{equation}
\left|\frac{1}{2i} \int_{\gamma_2} f \, \dtext \bar z \right| \leqslant \frac{1}{2} \int_{|z|=\epsilon} |f|\, |\dtext z| = \pi \epsilon \| f \|_{\infty}
\end{equation}
Passing to the limit as $\epsilon \rightarrow 0$ we find that
\begin{equation}
\int_{\mathbb U}  \left(\re f_z \right)\, \dtext x \, \dtext y =0
\end{equation}
which gives the desired contradiction.

Therefore, there are no elliptic sectors in $\Omega$. We now come to the fundamental theorem of Brouwer \cite[Theorem 5]{Br} which asserts that the index of $f$ at the point $0$ is given  by
\begin{equation}
n_f(0) = 1 + \frac{n_e-n_h}{2}
\end{equation}
where $n_e$ stands for the number of elliptic sectors and $n_h\geqslant 0$ stands for the number of hyperbolic sectors in $\Omega$. We just proved that $n_e=0$. Since $n_h \geqslant 0$ and $n_f(0) \geqslant 1$, this is only possible if $n_h=0$ and $n_f(0)=1$, as claimed. \qed

\section{Proof of Theorem \ref{T15}}\label{PT15}
Let $f=u+iv$. In this notation the condition $\re f_z =0$ reads as $u_x+v_y=0$. Therefore there exists a real valued function $\psi$ such that
$${ \psi_x} = - v \hskip0.5cm   \textnormal{ and } \hskip0.5cm  {\psi_y} = u $$
or, using  complex notation,  $\nabla \psi =\psi_x +i\psi_y= i f$.  For $(a,b) \in \Omega $ we define
$$\psi^{a,b}(x,y)= \psi(x,y)- \big[ \psi (a,b) +(x-a) \psi_x (a,b) + (y-b) \psi_y (a,b)\big] $$
Due to the local injectivity of $f$, $(a,b)$ is an isolated critical point of $\psi^{a,b}$.
Since the topological index of $\nabla \psi^{a,b} (a,b)$ is equal to $1$, by \cite[Lemma 3.1]{AM} there is a neighborhood $U$ of $(a,b)$ such that either
\begin{enumerate}[(i)]
\item\label{above} $\psi^{a,b}>0$ in $U \setminus (a,b)$, or
\item\label{below} $\psi^{a,b}<0$ in $U \setminus (a,b)$.
\end{enumerate}
We claim that only one of the above alternatives occurs for all $(a,b)\in \Omega$. Suppose to the contrary that
\eqref{above} occurs at $(a_1,b_1)$ and \eqref{below} occurs at $(a_2,b_2)$. Consider the function
$\phi(t)=\psi(a_1+t(a_2-a_1),b_1+t(b_2-b_1))$ which is defined on some open interval containing $[0,1]$, because $\Omega$ is convex.
Since any tangent line to the graph of $\phi$ stays (locally) on one side of the graph,
the Mean Value Theorem implies that $\phi'$ does not have any points of local extremum. Therefore, $\phi'$ is monotone,
and $\phi$ is either convex or concave. However, this contradicts our assumptions about $(a_1,b_1)$ and $(a_2,b_2)$.

Suppose for the sake of definiteness that only~\eqref{above} occurs in all domain $\Omega$. It follows that $\psi$ is strictly
convex in $\Omega$. Being the gradient mapping of a strictly convex function, the map $if$ is injective~\cite[Corollary 26.3.1]{Rob}, and so is $f$. \qed

\section{Proof of Theorem \ref{T1}}
Once we know that $f$ is locally quasiconformal, by Theorem~\ref{T2}, its global injectivity  is a consequence of integral estimates near $\infty$.  The following elementary, though interesting fact yields Theorem~\ref{T1}.
\begin{proposition}\label{P1}
If $f :\mathbb C \rightarrow \mathbb C$ is locally $K$-quasiconformal and $J_f(z)=\left|f_z\right|^2 - \left|f_{\bar z}\right|^2 \geqslant \lambda^2 $, almost everywhere for some $\lambda >0$, then $f$ is injective. Precisely we have
\begin{equation}\label{Eq31}
\left|f(z_1)- f(z_2)\right| \geqslant \frac{\lambda}{\sqrt{K}}\,  |z_1 -z_2|
\end{equation}
\end{proposition}
\begin{proof} We may assume that $f(0)=0$ and $f(1)=1$, by rescaling if necessary. Stoilow factorization provides us with a normalized  $K$-quasiconformal map $\chi : \mathbb C \rightarrow \mathbb C$, $\chi (0)=0$, $\chi (1)=1$, such that $H(\omega)= f \big(\chi (\omega)\big)$ is an entire function. We aim to show that $H(\omega) \equiv \omega$. Since $f$ is locally injective so is $H$. In particular, $H^\prime (\omega) \not= 0$. By the chain rule we have the following lower bound of the derivative
$$\left|H^\prime (\omega)\right|^2 = J_f(z) \, J_\chi (\omega) \geqslant \lambda^2 \, J_\chi (\omega)$$
Choose and fix a sufficiently small positive number $\epsilon$, for instance $0< \epsilon < \frac{1}{K-1}$ will suffice. Then we have
$$\frac{1}{|H^\prime (\omega)|^{2 \epsilon}} \leqslant \frac{1}{\lambda^{2 \epsilon} \, J_\chi^\epsilon (\omega)}$$
Consider the entire function
$$F(\omega)= \left[H^\prime (\omega)\right]^{-\epsilon} = \sum_{m=0}^\infty a_m \, \omega^m$$
Integration over the disk $B= \{\omega \, ; \; \; |\omega| \leqslant R\}$, yields
$$\sum_{m=0}^\infty \frac{|a_m|^2}{m+1}R^{2m}= \dashint_B \left| F(\omega)\right|^2 \, \dtext \omega \leqslant \frac{1}{\lambda^{2\epsilon}}  \dashint_B \frac{\dtext \omega}{J^\epsilon_\chi (\omega)}$$
The integral average in the right hand side exhibits a power growth with respect to $R$. Although the precise value of the power is immaterial for the forthcoming arguments, we demonstrate here the use of Astala's area distortion theorem, \cite{As}, to obtain sharp power.
\begin{lemma}
If $\chi : \mathbb C \rightarrow \mathbb C$ is $K$-quasiconformal and $B \subset \mathbb C$ is a disk, then for every $0< \epsilon < \frac{1}{K-1}$ we have
$$\dashint_B \frac{\dtext \omega}{J_\chi^\epsilon (\omega)} \leqslant \frac{C_K}{1-(K-1)\epsilon } \left(\frac{|B|}{|\chi (B)|}\right)^\epsilon$$
In particular, if $\chi(0)=0$ and $\chi(1)=1$, then for $R \geqslant 1$
$$\frac{1}{\pi \, R^2}\int_{|\omega| \leqslant R} \frac{\dtext \omega}{J^\epsilon_\chi (\omega)} \leqslant C_K(\epsilon) R^{2 \, \epsilon \, (1-1/K)}$$
\end{lemma}
\begin{proof} See \cite[Theorem 13.2.4]{AIMbook}.
\end{proof}

We just arrived at the inequality
$$\sum_{m=0}^\infty \frac{|a_m|^2}{m+1}R^{2m} \leqslant \frac{C_K(\epsilon)}{\lambda^{2\epsilon}} R^{2 \, \epsilon \, (1-1/K)}$$
where $2 \, \epsilon \, (1-1/K) < 2/K <2$. Therefore, $a_m =0$ for $m \geqslant 1$. This means that $F(\omega)$ is constant, and  so is $H^\prime (\omega)$. Hence $H(\omega)= \omega$, because of normalization $H(0)=0$ and $H(1)=1$. In conclusion, $f(z)$ is inverse of $\chi (\omega)$, and we have
$$\left|D \chi (\omega)\right|^2 \leqslant K J_\chi (\omega) = \frac{K}{J_f (z)} \leqslant \frac{K}{\lambda^2}$$
$$\left| \chi (\omega_1) - \chi (\omega_2)\right| \leqslant \frac{\sqrt{K}}{\lambda} |\omega_1 - \omega_2|$$
which is equivalent to (\ref{Eq31}).
\end{proof}
Returning to Theorem \ref{T1}, we consider quasiregular mappings
$$f^\lambda (z)= f(z)+ \lambda z\, , \hskip1cm \lambda \geqslant 0$$
Clearly,
$$J_{f^\lambda} (z)= J_f(z)+ \lambda^2 +2\lambda \re f_z \geqslant \lambda ^2$$
Hence $f^\lambda$ is $K$-quasiconformal mapping of $\mathbb C$ onto itself, for all $\lambda >0$. Passing to the limit as $\lambda \rightarrow 0$, by  Theorem \ref{HuT} we conclude that $f$ is injective in the entire plane.\qed

\section{Examples}

\begin{example}\label{branchexample}
For every $\epsilon >0$ there is a nonconstant quasiregular map $f\colon \mathbb C\to\mathbb C$ whose $z$-derivative lies in the sector
\begin{equation}\label{Ex1f1}
\re f_z \geqslant - \epsilon \left| \im f_z \right| \hskip1cm \mbox{ a.e.  in } \mathbb C
\end{equation}
and yet $f$ fails to be injective.
 \end{example}
\begin{proof}
We need only consider $0<\epsilon \leqslant 2$. Let us introduce a parameter $\delta = \frac{\epsilon}{ 2 + \sqrt{4-\epsilon^2}} \leqslant \frac{\epsilon}{2} \leqslant 1$ so that $\epsilon = \frac{4 \delta }{1+\delta^2}$. First we define a quasiconformal homeomorphism of the upper half plane $\mathbb H=\{z\in\mathbb C\colon \im z> 0\}$ onto the complex plane with a slit along the nonnegative $x$-axis.
\begin{equation}
f(z)= \begin{cases}  \displaystyle \frac{2z^2}{|z| \sqrt{1+ \delta^2}},\quad &\textnormal{if } \re z \geqslant - \delta \im z \\ (i-\epsilon )z - i \bar z,\quad &  \textnormal{if } \re z \leqslant - \delta \im z
  \end{cases}
\end{equation}
A straightforward computation shows that
$$(i- \epsilon) z - i \bar z =  \frac{2 z^2}{|z| \sqrt{1+\delta^2}} \hskip0.5cm \textnormal{on the ray } z=(i-\delta)t \, , \; \; t>0$$
Thus $f$ is continuous on $\mathbb H$. Moreover, its complex derivatives outside this ray are:
$$f_z = \begin{cases}\frac{3z}{|z| \sqrt{1+\delta^2}},\quad &\textnormal{if } \re z > - \delta \im z  \\
i-\epsilon,\quad &\textnormal{if } \re z <- \delta \im z       \end{cases}$$
$$f_{\bar z} = \begin{cases}\frac{-z^3}{|z|^3 \sqrt{1+\delta^2}},\quad &\textnormal{if } \re z > - \delta \im z  \\
-i,\quad &\textnormal{if } \re z < - \delta \im z       \end{cases}$$
In any case we find that
\begin{equation}\label{Exine3}
\left|f_{\bar z} \right| \le \frac{1}{\sqrt{1+\epsilon^2}} \left|f_{z} \right|
\end{equation}
Regarding the condition (\ref{Ex1f1}), we have
$$\re f_z \ge \begin{cases}-\delta \left| \im f_z\right|,\quad &\textnormal{if } \re z > - \delta \im z  \\
-\epsilon  \left| \im f_z \right|,\quad &\textnormal{if } \re z < - \delta \im z       \end{cases}$$
In any case, $\re f_z \ge - \epsilon \left| \im f_z \right| $. Next we note  that $f\colon \mathbb H \to\mathbb C$ extends continuously to the real line with values  in   the nonnegative real axis. Precisely, we have
$$f(x+i\, 0) = \begin{cases}\frac{2x}{\sqrt{1+\delta^2}},\quad &\textnormal{for } x \ge 0  \\
-\epsilon x,\quad &\textnormal{for }  x \le 0      \end{cases}$$
By  reflection about the $x$-axis, we extend $f$ to the lower half-plane; that is, by setting $f(z) = \overline{f(\bar z)}$ if $\im z \le 0$. Inequality (\ref{Exine3}) holds almost everywhere in $\mathbb C$. But $f$ is no longer injective near any neighborhood of the origin, because
$$f(z)=f(-z) \hskip0.5cm \textnormal{ if } \re z \ge - \delta |\im z|$$
It only remains to verify (\ref{Ex1f1}) in  the lower half plane. For $\im z <0$ we have $f_z(z)= \overline{f_z(\bar z)}$. Hence
$$\re f_z (z)= \re f_z(\bar z) \ge - \epsilon\,  \left|  \im f_z (\bar z) \right|  = -\epsilon \left| - \im f_z (z)\right| = -\epsilon \left|  \im f_z (z)\right|$$
as desired.
\end{proof}
\begin{example}\label{noninjexample}
Let $\Omega=\{z\in\mathbb C\colon \re z> 0\}$ be the right halfplane. For every $M \geqslant 1$ there is a nonconstant quasiregular map $f\colon \Omega \to\mathbb C$
such that
\begin{equation}\label{Exf2}
\re f_z \geqslant M \left| \im f_z \right|
\end{equation}
a.e. in $\Omega$ and $f$ is not injective.
\end{example}
\begin{proof}
For a given $M \geqslant 1$ we define
\[
f(z)=\begin{cases}
(4M^2 +4Mi)z + (4M^2 +1) \bar z,\quad &  \textnormal{if } 0 < \re z \leqslant 2M \im z  \\
(4M^2 -4Mi)z + (4M^2 +1) \bar z     ,\quad & \textnormal{if } 0 < \re z \leqslant -2M \im z \\
(8M^2-1)z, \quad &  \textnormal{if } \re z \geqslant 2M \left| \im z\right| . \\
\end{cases}
\]
The reader may wish to verify that $f$ is  quasiregular and  satisfies inequality (\ref{Exf2}). Moreover, $f$ fails to be injective, because
$f(1+4Mi)=1-8M^2 =f(1-4Mi)$.
\end{proof}

\begin{example}\label{gradexample} For any $n\ge 3$ there exists a function $\psi\in C^1(\mathbb R^n)$
such that $\nabla \psi\colon \mathbb R^n\to\mathbb R^n$ is a nonconstant quasiregular mapping that is not a local
homeomorphism.
\end{example}
\begin{proof} Define
\[\psi(x_1,\dots,x_n)=\frac{(x_1^2-x_2^2)^2-4x_1^2x_2^2}{x_1^2+x_2^2}-\frac{x_3^2}{2}+\frac12\sum_{k=4}^n x_k^2.\]
Note that $\psi$ is homogeneous of degree $2$. Since $\psi$ is smooth away from the origin, the homogeneity
implies that the entries of the Hessian matrix $D^2\psi$ belong to $L^{\infty}(\mathbb R^n)$. We claim
that $\det D^2\psi\ge 16$ a.e. Thanks to the block-diagonal form of $D^2 \psi$ and $\prod_{k=3}^n \frac{\partial^2 \psi}{\partial x_k^2}=-1 $ it suffices to show that
$\det D^2 u\le -16$, where $u$ is the restriction of $\psi$ to the plane $\{x\colon x_k=0 \ \forall k\ge 3\}$.
Writing $u$ in complex notation,
\[u(z)=\frac{\re(z^4)}{\abs{z}^2}=\frac{1}{2}\left(\frac{z^3}{\bar z}+\frac{\bar z^3}{z}\right)\]
we find that
\begin{equation}\label{uzbar}
\nabla u=2u_{\bar z}=\frac{3\bar z^2}{z}-\frac{z^3}{\bar z^2}
=\frac{\bar z^3}{\abs{z}^2}\left(3-\frac{z^4}{\bar z^4} \right)
\end{equation}
and
\[\det D^2u = 4(\abs{u_{\bar z z}}^2-\abs{u_{\bar z \bar z}}^2)=-22-6\re\frac{z^4}{\bar z^4}\le -16.\]
This proves that $\nabla \psi$ is quasiregular and, moreover, belongs to the class BLD (bounded length distortion) introduced by Martio and V\"ais\"al\"a~\cite{MV}. The last part of~\eqref{uzbar} shows that $\nabla u$ is homotopic to $z\mapsto \bar z^3$
in $\mathbb C\setminus\{0\}$. Since $\nabla u$ has index $-3$ at $0$, it follows that $\nabla \psi$ has index $3$
at all points $x\in\mathbb R^n$ with $x_1=x_2=0$. In fact, these points constitute the brach set of $\nabla \psi$, because $\nabla \psi$ is a local diffeomorphism.
\end{proof}

\section*{Acknowledgement}
We would like to  thank the referee  for a careful reading of the manuscript and several  valuable comments which  helped us to improve the presentation.

\bibliographystyle{amsplain}

\end{document}